\documentclass{article}
\usepackage[centertags]{amsmath}
\usepackage{hyperref}
\usepackage{amsfonts}
\usepackage{amssymb}
\usepackage{amsthm}
\usepackage{newlfont}
\usepackage{amscd}
\usepackage{amsmath,amscd}
\usepackage[all,knot,arc,poly]{xy}
\usepackage{verbatim}
\usepackage{enumerate}
\usepackage{color}
\input xy

\xyoption{all}
\usepackage[neverdecrease]{paralist}

\newcommand{\A}{\mathbb{A}}

\def\Om{{\Omega}}

\newtheorem{thm}{Theorem}[section]
\newtheorem{cor}[thm]{Corollary}
\newtheorem{lem}[thm]{Lemma}
\newtheorem{prop}[thm]{Proposition}
\newtheorem{pro}[thm]{Proposition}
\newtheorem*{thma}{Theorem A}
\newtheorem*{thmb}{Theorem B}
\newtheorem*{thmc}{Theorem C}

\theoremstyle{definition}

\newtheorem{example}{Example}
\newtheorem{defn}[thm]{Definition}

\newtheorem{notation}[thm]{Notation}
\newtheorem{observation}[thm]{Observation}

\theoremstyle{remark}
\newtheorem{rem}[thm]{Remark}
\newtheorem{remarks}[thm]{Remarks\newline\newline}

\DeclareMathOperator{\holim}{holim}

\def\rar{\rightarrow}
\def\lar{\leftarrow}
\def\lrar{\longrightarrow}

\def\hrar{\hookrightarrow}

\def\drar{\dashrightarrow}
\def\bu{\bullet}

\newcommand\ackname{Acknowledgements:}
\if@titlepage
  \newenvironment{acknowledgements}{%
      \titlepage
      \null\vfil
      \@beginparpenalty\@lowpenalty
      \begin{center}%
        \bfseries \ackname
        \@endparpenalty\@M
      \end{center}}%
     {\par\vfil\null\endtitlepage}
\else
  
\fi

\title{Homotopy Normal Maps}

\author{Matan Prezma\footnote{later papers appear under ``Matan Prasma"}}
\date{}

\begin{document}

\maketitle

\begin{abstract}
A group property made homotopical is a property of the corresponding classifying space. This train of thought can lead to a homotopical definition of 
normal maps between topological groups (or loop spaces). 

In this paper we deal with such maps, called \emph{homotopy normal maps}, which are topological group maps $N\rar G$ being `normal' in that they induce a compatible topological group structure on the homotopy quotient $G//N:=EN\times_N G$. We develop the notion of homotopy normality and its basic properties, and show it is invariant under homotopy monoidal endofunctors of topological spaces, e.g. localizations and completions. In the course of characterizing normality, we define a notion of a \emph{homotopy action} of a loop space on a space phrased in terms of Segal's 1-fold delooping machine. Homotopy actions are `flexible' in the sense they are invariant under homotopy monoidal functors, but can also rigidify to (strict) group actions.  
\end{abstract}

\begin{ackname}
This paper is based on the author's Ph.D. thesis at the Hebrew University of Jerusalem.
The author would like to express deep thanks to his advisor, Emmanuel Dror Farjoun for his continuous guidance, discussions and encouragement. The author would also like to thank the Hebrew University of Jerusalem for support of his studies. Special thanks are extended to David Blanc and James Stasheff for helpful suggestions and conversations.
\end{ackname}

\section{Introduction}
Homotopy normality is an attempt to derive a homotopical analogue for the inclusion of a normal subgroup via classifying spaces.
An inclusion of topological groups $N\hrar G$ is the inclusion of a normal subgroup if and only if it is the 
kernel inclusion of some group map $G\rar H$. Since any map is, up to homotopy, 
an inclusion, one needs to consider all group maps $N\rar G$. Such a map should then be `homotopy normal' if $BN\rar BG$ is 
the map from the homotopy fiber to the total space for some map $BG\drar W$. There is another angle from which this 
notion makes sense. To every group map $N\rar G,$ one can associate the Borel construction 
$EN\times_N G=:G//N$, which is the `correct' quotient in the homotopical world.
We note that such an extension $BG\drar W$ induces a loop space structure on $G//N$, and a loop map 
structure (up to map equivalence) on $G\rar G//N$, providing a second analogy to the group theoretic 
notion: a group inclusion $N\hrar G$ is the inclusion of a normal subgroup if and only if $G/N$ admits a group structure for which
the natural quotient map $G\rar G/N$ is a group map.

Let $f:X\rar Y$ be a pointed map of connected spaces. 
Consider the Puppe-Nomura sequence \cite{Nom}
$$
\Omega X\rar \Omega Y \rar \Omega Y//\Omega X\rar X\rar Y,
$$ 
where we denote $\Omega Y//\Omega X:=hfib(f)$.

The following is essentially taken from \cite[$\S$5]{FS}.
\begin{defn}\label{homotopy normal}
A loop map $\Omega f:\Omega X{\rar} \Omega Y$ is \emph{homotopy normal} if 
there exist a connected space $W$ with a map 
$\pi :Y\rar W$, so that $\xymatrix{X\ar[r]^{f} & Y\ar[r]^{\pi} & W}$ 
is a homotopy fibration sequence.
The map $\pi:Y\rar W$ is called a \emph{normal structure}.
\end{defn}

\begin{remarks}\label{remarks}
\begin{enumerate}[(a)]
\item We see that a loop map $\Omega f:\Omega X\rar \Omega Y$ is homotopy normal if and 
only if $f:X\rar Y$ admits a structure of a homotopy principal fibration, 
i.e. equivalent to a principal fibration. In particular, the homotopy fiber of such a loop map has the structure of a double loop space. 

\item If $\Omega f:\Omega X\rar \Omega Y$ is homotopy normal, the group map\linebreak
$\pi_0(\Om f):\pi_0(\Om X)\rar \pi_0(\Om Y)$ is normal in the sense of \cite{FS} i.e. underlies a crossed module structure on the corresponding groups. Whitehead showed (see \cite{WH}) that crossed modules correspond to connected 2-types. We note that if a discrete group map $N\rar G$ is normal (in the sense of \cite{FS}) and $BG\rar W$ its normal structure then $W$ is the corresponding connected 2-type. 

\end{enumerate}

\end{remarks}

\begin{example}
If $F\rar E\rar B$ is a fibration sequence, the map $\pi_1F\rar\pi_1E$ is a homotopy normal map of discrete groups. It is also true that any homotopy normal map of discrete groups is of this form (see \cite[\S 2.6]{BHS} and \cite[corollary 1.5]{Lod}).
\end{example}

\begin{example}\label{double loop map}
Any double loop map $\Omega^2 f:\Omega^2 X{\lrar} \Omega^2 Y$ where $X,Y$ are simply connected spaces 
is homotopy normal: take $W=hfib(X\rar Y) $; 
$W$ is then a connected space which extends the Puppe-Nomura sequence.
\end{example}

\begin{example}
Let $F$ be a pointed connected space. Then the universal fibration in \cite{Got}, $F\rar Baut_*(F)\rar Baut(F)$ induces a homotopy normal map $\Om F\rar \Om B aut_*(F)$. This map may be viewed as a universally initial homotopy normal map in the following sense:
for every homotopy normal map $\Om F\rar \Om X$ there exist a loop map $\Om X\rar \Om Baut_*(F)$ and a homotopy commutative triangle
$$\xymatrix{\Om F\ar[r]\ar[rd] & \Om Baut_*(F).\\ & \Om X\ar@{-->}[u]}$$
The dashed arrow is obtained as follows. Assume $F\rar X\rar W$ is a homotopy fibration sequence giving a normal structure on $\Om F\rar \Om X$. By \cite{Got}, there exist a `classifying map' $c:W\rar Baut(F)$, such that $X\rar W$ is obtained as a homotopy pullback $$\xymatrix{X\ar[r]\ar[d] & W\ar[d]^c\\Baut_*(F)\ar[r] & Baut(F).}$$  
This can be extended to a homotopy commutative diagram $$\xymatrix{F\ar[d]^{\simeq}\ar[r] & X\ar[r]\ar[d] & W\ar[d]^c\\F\ar[r] & Baut_*(F)\ar[r] & Baut(F)}$$
and looping down $X\rar Baut_*(F)$ gives the desired map.
\end{example}

\subsection*{Main results}
Given a group map $N\rar G$, each level of the bar construction 
$Bar_\bullet (G,N)=\{G\times N^k\}_{k\geq 0}$ \cite[$\S$7]{May} admits an action of $G$, namely
the one induced from the group inclusions $s_0:G\rar G\times N$, $s_1s_0:G\rar G\times N^2$, etc.
Similarly, in any simplicial group $\Gamma_\bullet$, $\Gamma_0$ 
acts on each level via degeneracies (as above),
and endows $\Gamma_\bullet$ with a structure of $\Gamma_0-$simplicial set. 

The following is the main theorem in \cite[$\S$4]{FS}, rephrased.
\begin{thm} \label{[FS]}
A map of discrete groups $f:N\rar G$ is homotopy normal if and only if there exists a simplicial 
group $\Gamma_\bullet$, with an isomorphism $\Gamma_0\cong G$ which extends to a $G$-equivariant 
isomorphism of simplicial sets 
$$
Bar_\bullet(G,N)\rar \Gamma_\bullet.
$$
\end{thm}

The main goal of this work is to describe a generalization of theorem \ref{[FS]} that characterizes
all normal maps $\Om X\rar \Om Y$. Our strategy is as follows.

In $\S$\ref{fiber} we define a homotopical analogue to the bar construction in the case of loop maps 
$\Omega X\rar \Omega Y$, $Bar_\bullet(\Omega Y,\Omega X)$. In the degenerate case of 
$\Omega Y \simeq *$, $Bar_\bullet(*,\Omega X)=Bar_\bullet(\Omega X)$, and 
one recovers Segal's 1-fold delooping machine (Definition \ref{special}) for $\Omega X$.

Next, in \S\ref{homotopy action section} we define the notion of a homotopy action of a loop space on a space. We study its
basic properties, and establish a weak equivalence between the category of homotopy actions of a fixed loop space and the category of spaces with an action of a fixed topological group. The simplicial space $Bar_\bu(\Om Y, \Om X)$ admits a canonical homotopy action of $\Om Y$. A homotopy action of $\Om Y$ is also defined for any simplicial loop space $\Gamma_\bu$ satisfying $\Gamma_0\simeq \Om Y$ . Using this setup we can state a homotopical analogue of theorem \ref{[FS]}.

\begin{thma}\label{characterization of normality}
A loop map $\Omega f:\Omega X\rar \Omega Y$ is homotopy normal 
if and only if there exists a simplicial loop space $\Gamma_\bullet$ with $\Gamma_0 \simeq \Omega Y$ (as loop spaces), and such 
that the canonical homotopy actions of $\Omega Y$ on $\Gamma_\bullet$ and on
$Bar_\bullet(\Omega Y,\Omega X)$ are weakly equivalent.
\end{thma}

As often happens, theorem \ref{[FS]} is a special case of theorem A in that it is precisely its $\pi_0$ statement.
One consequence of theorem A is the fact that homotopy normal maps are invariant under homotopy monoidal functors.

\begin{defn}
A functor $L: Top\rightarrow Top$ is called a \emph{homotopy monoidal} (HM) functor if it preserves 
homotopy equivalences, contractible spaces, and finite products up to homotopy. The last condition can also be formulated as follows:
for every pair of spaces $X,Y$, the canonical map $L(X\times Y)\overset{\simeq}{\rar} LX\times LY$ is a homotopy equivalence.
\end{defn}

Let L be an HM functor and $\Omega f: \Omega X\rar \Omega Y$ a loop map.
It is implicit in \cite{Bous} and \cite{Far}, and can be proved also by using the delooping theorem of \cite{Seg} 
that $L(\Omega X)$ always has the homotopy type of a loop space and $L(\Omega f)$ is always equivalent to a loop map.

\begin{rem}
Although HM functors preserve the property of having (the homotopy type of) a loop space, they do not commute with the functor \linebreak $\Om:Top_*\rar Top_*$. 
\end{rem}

Using the fact that homotopy actions of loop spaces, can be described in terms of maps between finite products of spaces we show that HM functors preserve homotopy normality.

\begin{thmb}
Let $\Omega f:\Omega X\!\rar\Omega Y$ be a homotopy normal map. 
If \mbox{$L\!:\!Top\!\rightarrow\! Top$} is an HM functor, then 
$L(\Om f):L\Omega X\rar L\Omega Y$ 
is a homotopy normal map.
\end{thmb}

This, in turn, gives an immediate proof of a theorem due to Dwyer and Farjoun (\cite[$\S$3]{DF}) which we restate.

\begin{thmc}
Let $f:X\rar Y$ be a map of pointed connected spaces and 
$p:E\rar B$ be a homotopy principal fibration of connected spaces. 
If $L_{\Sigma f}$ is the localization functor by $\Sigma f :\Sigma X\rar \Sigma Y$, 
then ${L_{\Sigma f}(p)}:L_{\Sigma f}E{\lrar} L_{\Sigma f} B$ 
is a homotopy principal fibration.
\end{thmc}

\begin{rem}
In what follows, we use $L$ to denote an arbitrary HM functor. The notation $L$ reflects the special 
case of localization by a map.
\end{rem}

We would like to refer the reader to related work, \cite{FH}, on homotopy\linebreak (co-)normal structures in a category with a class of weak equivalences and some additional structure, called a twisted homotopical category.


\section{Preliminaries}
Throughout this paper, \emph{topological spaces} or \emph{spaces} will mean topological spaces of the homotopy type of CW complexes. We denote the corresponding category by $Top$. Thus, by Whitehead's theorem, every weak equivalence is in fact a homotopy equivalence.
All mapping spaces will be taken with the compact-open topology.
The \emph{path space} $PX$ of a pointed space $X$ is the space of maps $\{\alpha:I\rar X | \alpha(0)=*\}$;
a \emph{loop space} is understood to be a space of the form $\Omega X:=\{\alpha:I\rar X | \alpha(0)=*=\alpha(1)\}$ where $X$ is a pointed connected space and a \emph{loop map} is a map of the form $\Om f:\Om X\rar \Om Y$ where $f:X\rar Y$ is a pointed map.
The following is a well-known fact, essentially contained in \cite{Kan} and \cite{Mil}.

\begin{thm}\label{rigidification}
If $X$ is a (pointed) connected space, then there exists a topological group $G$, 
with $X\overset{\simeq}{\rar} BG$. Moreover, one can construct $G$ functorially in $X$, i.e. 
if ${\Omega f}:\Omega X {\rar} \Omega Y$ is a loop map, there is a commutative diagram
 $$
\xymatrix{\Omega X \ar[r] \ar[d]_\simeq & \Omega Y\ar[d]_\simeq\\ G \ar[r] & H}
$$ 
with the vertical arrows being homotopy equivalences, 
and the bottom arrow being a topological group map. 
\end{thm}

A map $E\rar B$ is a (Serre) \emph{fibration} if it has the right lifting property with respect to all inclusions of the form $D^n\hrar D^n\times I$ that include the n-disc $D^n$ as $D^n\times\{0\}$.
A \emph{fibration sequence} is a sequence of the form $F\rar E\overset{p}{ \rar} B$, where $p:E\rar B$ is a fibration and either $(B,b_0)$ is pointed and
$F=p^{-1}(b_0)$ or $F=p^{-1}(b)$ for some $b\in B$ and $B$ is connected.
A sequence $X\rar Y\rar Z$ is called a \emph{homotopy fibration sequence} if there is a commutative 
diagram
\begin{equation*}
\begin{gathered}
\xymatrix{X\ar[r]\ar[d]_\simeq & Y\ar[r]\ar[d]_\simeq & Z\ar[d]_\simeq\\ F\ar[r] & E\ar[r] & B}
\end{gathered}
\end{equation*}
with vertical arrows being homotopy equivalences and the bottom being a fibration sequence. 
A homotopy fibration sequence $X\rar Y\rar Z$ is called a \emph{homotopy principal fibration sequence} if 
there is a connected space $B$ and a map $Z\rar B$, called the \emph{classifying map} such that $Y\rar Z\rar B$ is a homotopy fibration sequence. 
In that case, $X\simeq \Om B$ and there is a principal fibration sequence $G\rar E\rar E/G$, and a commutative diagram 
\begin{equation*}
\begin{gathered}
 \xymatrix{X\ar[r]\ar[d]_\simeq & Y\ar[r]\ar[d]_\simeq & Z\ar[d]_\simeq\\ G\ar[r] & E\ar[r] & E/G}
\end{gathered}
\end{equation*}
with all vertical maps being homotopy equivalences and the left vertical map being equivalent to a loop map $\Om B\rar \Om BG$.\\

As usual, we denote by $\Delta$ the category of finite ordinals
$[n]=(0,\ldots,n)$ with ordinal maps between them.
Given a category $\mathcal{C}$, a \emph{simplicial object} in $\mathcal{C}$ is a functor $\Delta^{op}\rar \mathcal{C}$, and we denote it by $X_\bullet$ with $X_n$ for its value on $[n]$.

Of special importance to this paper are simplicial objects in $Top$, namely 
\emph{simplicial spaces}. If $X$ is a space, we shall denote the constant simplicial space on it by $X$ when there is no risk of confusion.
An \emph{equivalence of simplicial spaces} (or: simplicial equivalence) is a simplicial map 
$f:X_\bullet\rar Y_\bullet$ such that, $f_n:X_n\rar Y_n$ is a homotopy 
equivalence for each $n$. Similarly, a \emph{(homotopy) fibration sequence of simplicial spaces} is a 
diagram of simplicial spaces $F_\bullet\rar E_\bullet\rar B_\bullet$ which is a level-wise 
(homotopy) fibration sequence.

We will often use a particular class of simplicial spaces introduced in a preprint of \cite{Seg} and originally called `group-like special $\Delta$-spaces'. Influenced by the terminology of \cite{Rez} we call them \emph{reduced Segal spaces}; these are defined as follows. 
\begin{defn}\label{special}(cf. \cite{Seg})

\begin{enumerate}[(a)]
\item
 A \emph{reduced Segal space} is a simplicial space $B_\bullet$ such that:
 \begin{enumerate}[(i)]
 	\item $B_0\simeq *$;
 	\item for each $n\geq 1$, the maps $p_n:B_n\rar B_1\times \cdots \times B_1$ (called Segal maps) induced by the maps 
 	$$
 	i_k:[1]\rar [n]\;\;\;\;\;(1\leq k\leq n)
 	$$ 
 	$$ 0\mapsto k-1\; ; \;1\mapsto k,$$
 	are homotopy equivalences;
 	\item the monoid structure on $\pi_0(B_1)$ admits inverses (i.e. is a group). 
 \end{enumerate} 
 
\item 
We say that $B_\bullet$ is a \emph{reduced Segal space for $\Om X$} if it comes equipped with a homotopy equivalence $|B_\bu|\overset{\simeq}{\rar} X$;
if  $B_\bu$ and $B_\bu'$ are reduced Segal spaces for $\Om X$, a map (respectively equivalence) between them is a simplicial map (respectively equivalence) $B_\bu\rar B_\bu'$ which makes the triangle of loop maps below commutative.
$$\xymatrix{\Om |B_\bu|\ar[rr]\ar[dr]_{\simeq} & & \Om |B_\bu'|\ar[dl]^{\simeq}\\ & \Om X & }$$ 
\end{enumerate}
\end{defn} 

\begin{rem}
By \cite[1.5]{Seg}, it follows that if $B_\bu$ is a reduced Segal space for $\Om X$ there is a natural homotopy equivalence $B_1 \overset{\simeq}{\rar} \Om |B_\bu|$. Thus, a reduced Segal space for $\Om Y$ can equivalently be defined as a reduced Segal space $B_\bu$ equipped with a loop equivalence $B_1\overset{\simeq}{\rar} \Om X$. The diagram of definition \ref{special}(b) should then be changed accordingly.
\end{rem}

For a topological group $G$ and $a:X\times G\rar X$ a right action of $G$ on a space $X$ which we denote by $x\mapsto xg$ for $x\in X$ and $g\in G$, the \emph{bar construction} (cf. \cite[\S 7]{May}) is the simplicial space $Bar_\bullet(X,G)$, consisting of:
\begin{enumerate}
\item
    for every $n\geq 0$, $Bar_n(X,G):=X\times G^n$ together with
\item
    face maps, $d_i^{(n)}\equiv d_i:Bar_n(X,G)\rar Bar_{n-1}(X,G)$ for every $n\geq 1$ and every $0\leq i\leq n$ given by:\\
    $d_i:(x,g_1,...,g_n)\mapsto \begin{cases}(x\cdot g_1,g_2,...,g_n) & if\;\; i=0\\ (x,g_1,..,g_{i-1},g_i\cdot g_{i+1},g_{i+2},...,g_n) & if \;\; 1\leq i <n \\ (x,g_1,...,g_{n-1}) & if \;\;i=n  \end{cases}$ \\
    and
\item
    degeneracy maps, $s_i:Bar_n(X,G)\rar Bar_{n+1}(X,G)$ for every $n\geq 1$ and every $0\leq i\leq n$ given by:\\
    $s_i:(x,g_1,...,g_n)\mapsto (x,g_1,...,g_i,e,g_{i+1},...,g_n)$.
\end{enumerate}

\section{The homotopy power of a map}\label{fiber}
Given a fibration $p:E\rar B$, one can define a simplicial space $Pow_\bullet(E\rar B)$, 
called the \emph{power of p}, by $Pow_n(E\rar B)=E\times_B E\cdots\times_B E \;\;\;(n+1\;times)$ with face 
and degeneracies being the obvious projections and diagonals. In \cite{Lod}, 
it is shown that for ($E$ non-empty and) $B$ connected, $|Pow_\bullet(E\rar B)|\simeq B$. We note that 
for a non-connected space $B$, $|Pow_\bullet(E\rar B)|$ is homotopy equivalent to the disjoint union of connected components of 
$B$ intersecting the image of $p$.

Here, we wish to construct such a power space for an arbitrary map
$f:X\rar B$ by means of homotopy pullbacks, thus turning it to a homotopically invariant construction.

We define the \emph{n-th homotopy power} of $f:X\rar B$ to be
$$hPow_n(X\rar B)=\xy<0.5cm,0cm>:
\POS (0,0) *+{map\Bigg(}
\POS (3,0) *+{ \xy <0.5cm,0cm>:
 	 \POS (0,2) *+{\Delta[n]_0} ="a"
 	 \POS (2,2) *+{X} ="b"
 	 \POS (1,1) *+{,}
 	 \POS (0,0) *+{\Delta[n]} ="c"
 	 \POS (2,0) *+{B} ="d"
 	 \POS "a" , \ar @{->}_\iota "c"
 	 \POS "b" , \ar @{->}^f "d"
 	 \endxy}
\POS (5,0) *+{\Bigg)} 
\endxy =\holim 
\xy<0.6cm,0cm>:
\POS (0,0) *+{\Bigg(}
\POS (2,0) *+{\xy<0.5cm,0cm>:
	\POS (0,2) *+{X} ="a"
	\POS (1,2) *+{X} ="b"
	\POS (2,2) *+{\cdots}
	\POS (3,2) *+{X} ="c"
	\POS (2,0) *+{B} ="d"
	\POS "a" , \ar@{->} "d"
	\POS "b" , \ar@{->} "d"
	\POS "c" , \ar@{->} "d"
\endxy}
\POS (4,0) *+{\Bigg)}
\endxy 	 ,
 	 $$                                  
with $\iota:\Delta[n]_0\rar\Delta[n]$ being the inclusion of the
$0$-skeleton into the topological $n$-simplex.

This clearly yields a functorial construction over $\Delta^{op}$, and we define:
\begin{defn}\label{hIFP}
The \textit{homotopy power} of a map $f:X\rar B$, denoted \linebreak
$hPow_\bullet(X\rar B)$, 
is the simplicial space with $hPow_n(X\rar B)$ on level $n$, and face and degeneracies 
given by the functorial construction above.
\end{defn}

Note that for a fibration $p:E\rar B$ one gets an equivalence of simplicial spaces $hPow_\bullet(E\rar B) \simeq Pow_\bullet(E\rar B)$. 

\begin{rem}
When calculating the homotopy power of a map $f:X\rar B$ we will often use a slightly different but equivalent construction. We first replace $f$ by an equivalent fibration $p:E_f\rar B$, i.e. one for which there is a commutative triangle $$\xymatrix{X\ar[d]_f\ar[r]^{\simeq} & E_f\ar[dl]^p \\ B}$$ and then take the power of $p$, as in \cite{Lod}. This construction is functorial as well. We also note that if $X\rar B$ is a pointed map, $hPow_\bu(X\rar B)$ naturally becomes a pointed simplicial space.
\end{rem}

\section{The homotopy bar construction}


Consider a topological group $G$ acting on a space $X$ and the corresponding (homotopy) principal 
fibration $G\rar X\rar X//G$. One has the `usual' bar construction 
$Bar_\bullet(X,G)=\{X\times G^k\}_{k\geq 0}$ with $|Bar_\bullet(X,G)|=X//G$. 
On the other hand, we can resolve $X//G$ by taking homotopy powers of the map \linebreak
$q:X\rar X//G$.  
\begin{pro}\label{Bar as IFP}
Let $G$ act on $X$ as above. Then there are simplicial equivalences 
$\xymatrix{
 Bar_\bullet(X,G)  \ar@<1ex>[r]
& hPow_\bullet(X\rar X//G) \ar@<1ex>[l] }$.
\end{pro}
\begin{proof}
Replacing $q:X\rar EG\times_G X$ by the fibration
$p:EG\times X\rar EG\times_G X$ and taking the pullback, we
get $hPow_1(X\rar X//G)=(EG\times X)\times_{X//G} (EG\times X)\cong EG\times G\times X$,
since $EG\times X$ is a free $G$-space. In general, 
$$
hPow_n(X\rar X//G)=(EG\times X)\times_{X//G} \cdots \times_{X//G}(EG\times X)
\cong EG\times X\times G^n,
$$
and the obvious map $EG\times X\times G^n \rar X\times G^n$ defines a
simplicial equivalence $hPow_\bullet(X\rar X//G) \rar Bar_\bullet(X,G)$. 
Taking (for example) Milnor's join construction, we have a natural base point for $EG$ and 
hence a canonical map $X\times G^n \rar EG\times X\times G^n$, which in turn defines another 
simplicial equivalence.
\end{proof}
In light of the last proposition, we define:
\begin{defn}\label{hBar}
Given a (homotopy) principal fibration sequence \mbox{$\Omega Y\!\rar \!X\!\overset{q}{\rar} \!Q$,} 
the \emph{homotopy bar construction} $Bar_\bullet(X,\Omega Y)$ is the homotopy power \linebreak
$hPow_\bullet(X\rar Q)$.   
\end{defn}
\begin{rem}\label{segal as special case}
In the case of a loop map $\Omega f: \Omega Y \rar \Omega Z$,  $Bar_\bullet(\Omega Z,\Omega Y)$ is 
the homotopy power of the map $q:\Omega Z \rar \Omega Z//\Omega Y:=hfib(f)$. If 
$\Omega Z\simeq *$, $Bar_\bullet(*,\Omega Y)$ becomes the power of the map $PY\rar Y$ which is 
a reduced Segal space for $\Omega Y$. Put differently, one can recover Segal's delooping machine by 
using homotopy powers.
\end{rem}
It is useful to have the following property.
\begin{pro}\label{loop of fiber product}
Let $f:X\rar B$ be any pointed map. The canonical map induces an equivalence
of simplicial spaces 
\mbox{$\Omega( hPow_\bullet(X\rar B)) \simeq$}   $hPow_\bullet(\Omega X\rar \Omega B)$.
\end{pro}
The proof is essentially the fact that given a pointed diagram 
$A\rar X \lar Y$, we have a weak equivalence $\Omega holim(A\rar X \lar Y)\simeq holim(\Omega A\rar \Omega X \lar \Omega Y)$.

\subsection{From homotopy normality to a simplicial loop space structure on the homotopy bar construction}
Let $\Omega f:\Omega X{\lrar} \Omega Y$ be a homotopy normal map.
We form the Puppe-Nomura sequence: 
$$
\xymatrix{\Omega X\ar[r]^{\Omega f} & \Omega Y\ar[r]^q &
\Omega Y//\Omega X\ar[r] &  X\ar[r] &  Y\ar[r]^\pi & W}.
$$ 
Then by \cite{Nom} there is a
commutative triangle in which the vertical arrow is a homotopy
equivalence 
\begin{equation*}
\begin{gathered}
\xymatrix{\Omega Y\ar[r]^q\ar[rd]_{\Omega \pi} & \Omega Y//\Omega X \ar[d]\\
                         & \Omega W}
\end{gathered}
\end{equation*}
Passing to (homotopy) powers, we get an equivalence
of simplicial spaces \linebreak
$$hPow_\bullet(\Omega Y\rar \Omega W)\simeq hPow_\bullet(\Omega Y\rar \Omega Y//\Omega X)$$
and, by proposition \ref{loop of fiber product}, an equivalence of
simplicial spaces 
$$\Omega (hPow_\bullet(Y\rar W))\simeq hPow_\bullet(\Omega Y\rar \Omega Y//\Omega X).$$
Using the argument above and definition \ref{hBar} we have just proved the following result.
\begin{thm}\label{Loop structure}
If $\Omega f:\Omega X\rar \Omega Y$ is homotopy normal, there are natural simplicial equivalences 
$\xymatrix{
 Bar_\bullet(\Omega Y,\Omega X)  \ar@<1ex>[r]
& \Omega(hPow_\bullet(Y\rar W) ) \ar@<1ex>[l] }$
\end{thm}
\begin{notation}\label{standard notation}(cf. \ref{Loop structure})

\begin{enumerate}
\item
For a homotopy normal map $\Omega f:\Omega X\rar \Omega Y$ and a given normal structure $\pi:Y\rar W$, we denote by $Q_\bullet$ the simplicial loop space

$\Omega(hPow_\bullet(Y\rar W) )$.
\item
The equivalences given in theorem \ref{Loop structure} will be denoted 
$$
\xymatrix{
 \epsilon:Bar_\bullet(\Omega Y,\Omega X)  \ar@<1ex>[r]
& Q_\bullet:\eta \ar@<1ex>[l] }.
$$
\end{enumerate}
\end{notation}
\begin{rem}
Notice that the maps $$\xymatrix{\epsilon_0:\Om Y  \ar@<1ex>[r]
& Q_0:\eta_0 \ar@<1ex>[l]}$$  are loop maps by construction, but for $n\geq 1$, the maps
$$\xymatrix{\epsilon_n:Bar_n(\Omega Y,\Omega X)  \ar@<1ex>[r]
& Q_n:\eta_n \ar@<1ex>[l]}$$
need not be loop maps. This means that we have, in general, two different loop space structures on $ \Om Y\times(\Om X)^n$. The non-trivial one is given by the equivalence $Bar_n(\Om Y,\Om X)\simeq Q_n$.
\end{rem}

\section{Homotopy actions}\label{homotopy action section}

By remark (a) in \ref{remarks} a homotopy normal map is a loop map with its underlying map having the 
structure of a principal fibration (of connected spaces). Furthermore, theorem \ref{[FS]} involves (strict) group actions. Hence, characterization and invariance of homotopy normal maps under HM functors should include characterization and invariance of group actions `up to homotopy' to some extent.   
Given an action of a topological group $G$ on a space $X$ and an HM functor $L:Top\rar Top$, we would 
like to construct a canonical `action' of $LG$ (not a group, not a loop space) on $LX$. In other 
words, we would like to have a homotopical notion of an action of (a space of the homotopy type of) a 
loop space on a space, invariant under HM functors. One approach we wish to refer the reader to is 
that of $A_\infty$-actions introduced in \cite{Now} and recently used in \cite{St11}. For our purpose, we could not use $A_\infty$-actions since it is not clear they are invariant under HM functors.  
As demonstrated in \S\ref{weakly correspondence}, homotopy actions can be rigidified into (strict) group actions. This rigidification gives in fact a `proxy action' on $X$ in the sense of \cite{DW} so all the homotopically-invariant information (e.g. homotopy fixed points) is preserved. Homotopy actions have more flexibility than proxy actions since the object which `acts' need not be a topological group but rather a loop space.
\subsection{Definition and basic properties}\label{construction}
If a topological group $G$ acts on a space $X$, one has a simplicial fibration sequence of the form
$X\rar Bar_\bullet (X,G)\rar B_\bullet G$, where the maps $X\rar Bar_n(X,G)$ and 
$Bar_n(X,G)\rar B_n G$ are given by $s_n\cdots s_0$ and projection respectively.

Under realization, this becomes a (homotopy) fibration sequence \linebreak
$X\rar X//G\rar BG$ with a connected base space, i.e. an `action up to homotopy' 
in the sense of \cite{DFK}.
The above simplicial fibration sequence is trivial in each level  $X\rar X\times G^n\rar G^n$, 
and hence constitutes a useful resolution. We note also that for all $n$, 
the map $d_1d_2\cdots d_{n}:Bar_n(X,G)\rightarrow Bar_0(X,G)$ is the projection on X and the map $d_0d_0\cdots d_0: Bar_n(X,G)\rar Bar_0(X,G)$ is given by $(x,g_1,...,g_n)\mapsto x\cdot (g_1\cdot ... \cdot g_n)$.

As we saw, the simplicial spaces $Bar_\bullet(X,G)$ and $B_\bullet G$ can be relaxed to their `homotopy versions', namely $Bar_\bullet(X,\Omega Y)$ and $Bar_\bullet(*,\Omega Y)$ (which is a reduced Segal space for $\Omega Y$ when $BG\simeq Y$).

\begin{defn}\label{defn h.action}
We say that a space $S$ of the homotopy type of a loop space, \emph{homotopy acts} on a space $X$, 
if there exist a simplicial map
$$
\xymatrix{ A_\bullet\ar[r]^\pi & B_\bullet}
$$ 
such that:
\begin{enumerate}
\item $A_0\simeq X$; \vspace*{4pt}
\item $B_\bullet$ is a reduced Segal space for $S$; 
\item for every $n$, the maps $\xymatrix{A_n\ar@{>}[rr]<1.2ex>^{d_1\cdots d_{n}\times \pi_n} \ar@{>}[rr]<-1.2ex>_{d_0\cdots d_0\times \pi_n}
 && A_0\times B_n}$ are 
homotopy equivalences.
\end{enumerate}
\end{defn}

Maps are defined as follows.
\begin{defn}
Given two homotopy actions of $S$ on $X$ and on $X'$, represented by $A_\bu\rar B_\bu$ and $A_{\bu}'\rar B_{\bu}'$ respectively, a map between them is a commutative square $$\xymatrix{A_\bu\ar[r]\ar[d] & B_\bu\ar[d]^{\simeq}\\ A_{\bu}' \ar[r] & B_{\bu}'}$$ such that the map $B_\bu\rar B_{\bu}'$ is an equivalence of reduced Segal spaces (see \ref{special}).
\end{defn}

\begin{notation}
We denote by ${Top}^{h\Om Y}$ the category of homotopy actions of (spaces of the homotopy type of) 
$\Om Y$ on spaces. 
\end{notation}
\begin{rem}\label{similar actions}
If $S\rar S'$ is a loop equivalence and $S$ homotopy acts on $X$, then $S'$ 
homotopy acts on $X$, since a reduced Segal space $B_\bu$ for $S$ induces a reduced Segal space for $S'$ simply by composing the map $B_1\overset{\simeq}{\rar} S$ with $S\overset{\simeq}{\rar} S'$ (see definition \ref{special}). 
\end{rem}

We will need a generalization of definition \ref{defn h.action} as follows.

\begin{defn}\label{homotopy action on a simplicial space}\label{equivalence}
A homotopy action of $\Omega Y$ on a simplicial space $X_\bullet$ is a
map of bisimplicial spaces $A_{\bullet\bullet}\rar B_{\bullet\bullet}$ 
such that for each $n$, $A_{\bullet n}\rar B_{\bullet n}$ 
is a homotopy action of $\Omega Y$ on $X_n$ and for every map $\theta:[n]\rar [m]$ in $\Delta$, $\theta^*: B_{\bu m} \rar B_{\bu n}$ is an equivalence of reduced Segal spaces for $\Om Y$; maps and equivalences are defined in the obvious way.
\end{defn}
 
\begin{observation}\label{Bar as haction}
If a topological group $G$ acts on a space $X$, the simplicial map 
$p:Bar_\bu(X,G)\rar B_\bu(G)$ is a homotopy action of $G$ on $X$. 
To see this, note that $B_\bu(G)$ is a reduced Segal space for $G$ and the maps 
$(d_1\cdots d_n)\times p_n:Bar_n(X,G)\rar Bar_0(X,G)\times B_n(G)$ 
are the identity maps $X\times G^n \rar X\times G^n$. One can verify that the maps 
$(d_0\cdots d_0)\times p_n:Bar_n(X,G)\rar Bar_0(X,G)\times B_n(G)$, i.e. the action of $G^n$ on $X$ (arising from multiplying $n$ elements in $G$ and then act on $X$) multiplied by the projection $p_n$, are homeomorphisms. 
\end{observation}

In \cite{Now}, the author defined an action of an $A_\infty$-space on a topological space. 
The difference between this approach and ours is essentially the difference between the approaches of 
Stasheff \cite{St70} and Segal \cite{Seg} to the characterization of loop spaces.

It is commonly said that in every fibration sequence, the loop space of the base `acts' on the fiber. 
We wish to demonstrate how a homotopy action interprets this statement.
\begin{thm}\label{classification of fibrations}
Given a fibration sequence $F\overset{i}{\rar} E\overset{p}{\rar} B$ with $B$ pointed connected, 
there is a homotopy action of $\Om B$ on $F$, 
represented by $A_\bullet\overset{\pi}{\rar} B_\bullet$, 
such that the map $|\pi|:|A_\bullet|\rar |B_\bullet|$ is equivalent to $p:E\rar B$.
\end{thm}
\begin{proof}
Consider the commutative square
\begin{equation*}
\begin{gathered}
\xymatrix{F\ar[r]\ar[d] & E\ar[d]\\ {*}\ar[r] & B}
\end{gathered}
\end{equation*}
Taking homotopy powers in each row produces a simplicial map 
$$
\pi:A_\bullet:= hPow_\bullet(F\rar E)\rar hPow_\bullet(*\rar B)=: B_\bullet.
$$ 
By remark \ref{segal as special case}, $B_\bullet$ is a reduced Segal space 
and thus $|B_\bullet|\simeq B$. Since $B$ is connected,  
it follows from $\S$\ref{fiber} that $|A_\bullet|\simeq E$.  
To see that $\pi:A_\bullet\rar B_\bullet$ is a homotopy action, we first replace 
$i:F\rar E$ and $*\rar B$ by equivalent fibrations \linebreak
$ev_1:F_i\rar E$ and $ev_1:PB\rar B$, 
where $PB$ is the path space and $F_i\subseteq F\times E^I$ is the space 
$\{(f,\alpha)|\alpha(0)=i(f)  \}$. 
Taking $\pi_0: F_i\rar PB$ to be 
$\pi_0(f,\alpha)=p\circ\alpha$ we obtain the commutative square 
\begin{equation*}
\begin{gathered}
\xymatrix{F_i\ar[r]^{ev_1}\ar[d]_{\pi_0} & E\ar[d]^p \\ PB\ar[r]^{ev_1} & B,} 
\end{gathered} \tag{$\ast$}
\end{equation*}
and taking powers (i.e. fiber products) of the rows, we obtain a simplicial map we denote as $\pi:A_\bullet \rar B_\bullet$.

Let us show that the maps

 $\xymatrix{A_1\ar@{>}[rr]<1.2ex>^{d_1\cdots d_{n}\times \pi_1} \ar@{>}[rr]<-1.2ex>_{d_0\cdots d_0\times \pi_1}
 && A_0\times B_1}$ are homotopy equivalences. We have a commutative cube

$$\xymatrix@!0{
A_1\ar[dd]_{d_1}\ar[rd]^{\pi_1}\ar[rr]^{d_0} & & F_i\ar'[d][dd]\ar[rd] \\
& B_1\ar[dd] \ar[rr] & & PB\ar[dd] \\ 
F_i\ar'[r][rr]\ar[rd] & & E\ar[rd] \\
& PB\ar[rr] & & B.
}$$
We want to show that the left-hand and upper faces are homotopy cartesian squares, which follows directly from the cartesianess of the lower, right-hand and outer faces using the fact that a square is cartesian if and only if the comparison map between homotopy fibers of rows/columns is a homotopy equivalence \cite[1.18]{GW}. 

One proceeds similarly to show that the maps $(d_0...d_0)\times \pi_n$ and $(d_1...d_n)\times \pi_n$ $(n>1)$ are homotopy equivalence.
Thus, $\pi:A_\bullet \rar B_\bullet$ is a homotopy action.

Lastly, since the equivalences $|Pow_\bullet(F_i\!\rar \!E)|\!\simeq \!E$ and 
\mbox{$|Pow_\bullet(PB\!\rar \!B)|\!\simeq \!B$} are natural, and in light of $(*)$ the map 
$|\pi|:|A_\bullet|\rar |B_\bullet|$ is equivalent to $p:E\rar B$. 
\end{proof}
The importance of theorem \ref{classification of fibrations} can be seen, for example, from the fact that it allows one to classify fibrations using homotopy actions.

Homotopy actions arise in our context in the following form. 
\begin{cor}\label{induced action}
If $\Om f:\Omega X \rar \Omega Y$ is a loop map, then $\Om f$ induces a homotopy action of $\Omega X$ on $\Omega Y$, natural in $f$.
\end{cor}

\begin{proof}
This follows from theorem \ref{classification of fibrations} if we consider the homotopy fibration sequence $\Om Y\rar \Om Y// \Om X \rar X$. 
Alternatively, if we (functorially) rigidify $\Om f:\Omega X \rar \Omega Y$ to a topological group map $G\rar H$ as in \ref{rigidification}, then as we saw, $Bar_\bullet(H,G)\rar B_\bullet G$ is a homotopy action. 
\end{proof}

Finally, let us see that homotopy actions are invariant under HM functors.
\begin{prop}
If $A_\bullet\rightarrow B_{\bullet}$ is a homotopy action of $\Omega Y$ on X, and \linebreak
$L:Top\rightarrow Top$ is an HM functor, then $LA_{\bullet}\rightarrow LB_{\bullet}$ 
 is a homotopy action of $L\Omega Y$ on $LX$.
\end{prop}
\begin{proof}
$LB_\bullet$ is a reduced Segal space for $LB_1$. 
In particular, $LB_1$ is of the homotopy type of a loop space.
Applying $L$ to the structure maps of the homotopy action yields the structure maps for 
$LA_\bullet\rar LB_\bullet$, and $L$ preserves homotopy equivalences.
\end{proof}
For the sake of completeness, we wish to define a map between homotopy actions of two non-homotopy-equivalent loop spaces. The simplicity of the definition demonstrates the `flexibility' of homotopy actions. For example, it allows one to talk about the category of \emph{all} homotopy actions.
\begin{defn}
 Given two homotopy actions of $\Omega Y$ on $X$ and of $\Omega (Y')$ on $X'$, represented by 
$\xymatrix{ A_\bullet\ar[r] & B_\bullet}$ and $\xymatrix{ A_{\bullet}' \ar[r] & B_{\bullet}'}$, 
a map between them is a commutative square of simplicial spaces 
\begin{equation*}
\begin{gathered}
\xymatrix{ A_\bullet\ar[r]\ar[d] & B_\bullet\ar[d] \\  A_{\bullet}'\ar[r] & B_{\bullet}'}
\end{gathered}.
\end{equation*}
Such a map will be called an \emph{equivalence} if both vertical maps are simplicial equivalences.
\end{defn}

\subsection{A weakly inverse correspondence with group actions}\label{weakly correspondence}
Our goal here is to establish a weakly inverse correspondence between the category $Top_{BG}$ of spaces over $BG$ and the category $Top^{h\Om Y}$ of homotopy actions of $\Om Y$ where $Y\simeq BG$. Since $Top_{BG}$ is Quillen equivalent to the category of $G$-spaces, we obtain a correspondence between homotopy actions and group actions which may be referred to as a `rigidification' of the homotopy action. 
Our functors will be weak inverses in the following sense.
\begin{defn}\label{weak equivalence}
Maps $f:X\rar Y$ and $f':X'\rar Y'$ are called \emph{weakly equivalent} if there is a zig-zag of commutative squares with all horizontal arrows being homotopy equivalences
\begin{equation*}
\vcenter{
\xymatrix{X\ar[d]_f\ar[r]^\simeq & X_1\ar[d] & \cdots \ar[l]_\simeq\ar[r]^\simeq & X_n\ar[d] & 
X'\ar[l]_\simeq\ar[d]^{f'}\\ 
Y\ar[r]^\simeq & Y_1 & \cdots\ar[l]_\simeq\ar[r]^\simeq & Y_n & Y'\ar[l]_\simeq}}
\end{equation*}
Similarly, simplicial maps $f:X_\bu\rar Y_\bu$ and $f':X_\bu'\rar Y_\bu'$ are 
called \emph{weakly equivalent} if there is a zig-zag of commutative squares as above, 
but with objects being simplicial spaces and maps being simplicial maps.
The number of squares involved in such a zig-zag is said to be its \emph{length}. In particular, maps are called \emph{equivalent} if they are weakly equivalent via a zig-zag of length $1$.
\end{defn}

\begin{defn}\label{functors}
Let $G$ be a topological group, $\Om Y$ a loop space and \mbox{$Y\rar BG$} a fixed homotopy equivalence.
\begin{enumerate} 
\item 
The functor $\mathcal{P}: Top_{BG}\rar Top^{h\Om Y}$ is defined as follows.
Given a map $E\rar BG$, let $X$ be its homotopy fiber. Thus, there is a commutative square $$\xymatrix{X\ar[r]\ar[d] & E\ar[d]\\PBG\ar[r]^{ev_1} & BG}$$
Then $\mathcal{P}(E\rar BG)$ is the map $hPow_\bu(X\rar E)\rar hPow_\bu(PBG\rar BG)$, which is a homotopy action of $\Om Y$ by proposition \ref{classification of fibrations}.
\item The functor $\mathcal{R}:Top^{h\Om Y}\rar Top_{BG}$ is defined as follows. Given a homotopy action $\pi: A_\bu\rar B_\bu$ of $\Om Y$ on $X$, $\mathcal{R}(A_\bu\rar B_\bu)$ is the composition $|A_\bu|\overset{|\pi|}{\rar} |B_\bu|\overset{\simeq}{\rar} Y\overset{\simeq}{\rar} BG$ where the second map comes from the fact that $B_\bu$ is a reduced Segal space for $\Om Y$(definition \ref{special}).

\end{enumerate}
\end{defn}
\begin{prop}\label{correspondence}
The functors above satisfy the following properties.
\begin{enumerate}[\upshape (a)]
\item If $E\rar BG$ is in ${Top}_{BG}$, then $\mathcal{P}(E\rar BG)$ is a homotopy action \linebreak of $\Om Y$ on $X:=hfib(E\rar BG)$.
\item If $\pi:A_\bu\rar B_\bu$ is a homotopy action of $\Om Y$ on $X$, then 
$\mathcal{R}(A_\bu\rar B_\bu)$ is a space over $BG$ with $X$ as its homotopy fiber.
\end{enumerate}
\end{prop}
\begin{proof}
\begin{enumerate}[\upshape (a)]
\item This follows from proposition \ref{classification of fibrations}.

\item Given a homotopy action $\pi:A_\bu\rar B_\bu$ of $\Om Y$ on $X$,
define a simplicial map $i:A_0\rar A_\bullet$ by
$i_n=s_{n-1}\cdots s_0$. Choose $b_0\in B_0$ and endow $B_n$ with a base-point 
$s_{n-1}\cdots s_0(b_0)$. By definition, the map 
$(d_1\cdots d_n)\times \pi_n:A_n\rar A_0\times B_n$ is a homotopy equivalence
and hence the map $\pi_n:A_n\rar B_n$ is equivalent to the trivial fibration
$A_0\times B_n\rar B_n$. We now claim that 
$A_0\overset{i_n}{\rar} A_n \overset{\pi_n}{\rar} B_n$ is a homotopy fibration sequence. 
To see this, note that by simplicial identities, the composite 
$\xymatrix{A_0\ar[rrr]^{(d_1 ... d_n\times \pi_n)\circ i_n} & & & A_0\times B_n}$ equals $1_{A_0}\times (\pi_n\circ i_n)$ and, since 
$B_0$ is contractible, $\pi_n\circ i_n=s_{n-1}\cdots s_0\circ \pi_0$ is null-homotopic. 
Hence, $i_n$ is equivalent to the fiber inclusion $A_0\rar A_0\times B_n$.
It follows that the sequence $A_0 \rar A_\bu \rar B_\bu$ is a 
homotopy fibration sequence in each level and so $A_0\rar |A_\bu|\rar |B_\bu|$ 
is a homotopy fibration sequence by \cite{Pup}. By definition, $A_0\simeq X$, and we are done.
\end{enumerate}
\end{proof}

\begin{thm}\label{homotopy action}
The functors $\mathcal{R}:\xymatrix{Top^{h\Om Y}\ar@{>}[rr]<1.2ex> && Top_{BG}\ar@{>}[ll]<1.2ex>}:\mathcal{P}$ of definition \ref{functors}
constitute a weakly inverse correspondence in the sense that:
\begin{enumerate}[\upshape (i)]
\item
$\mathcal{R}\mathcal{P}(E\rar BG)$ is weakly equivalent to $E\rar BG$;
\item
$\mathcal{P}\mathcal{R}(A_\bu\rar B_\bu)$ is weakly equivalent to $A_\bu\rar B_\bu$.
\end{enumerate}
\end{thm}

Theorem \ref{homotopy action} establishes a `rigidification theorem', which we wish to state separately.
\begin{thm}\label{bar rigidification}
Given a homotopy action of $\Om Y$ on $X$, represented by \linebreak
$\pi:A_\bu \rar B_\bu$, there is a topological group $G$ with $BG\simeq Y$ and a space $X'\simeq X$ together with a (strict) action of $G$ on $X'$ such that the simplicial map $\pi$ is weakly equivalent to the simplicial map $Bar_\bu(X',G)\rar B_\bu(G)$. 
\end{thm}

The proof of theorem \ref{homotopy action} will require some technical preparation. 
\begin{defn}
If $A_\bu$ is a simplicial space, the \emph{simplicial path space on $A_\bu$}, denoted $PA_\bu$, 
is the simplicial space defined by $PA_n=A_{n+1}$ with face maps $d_i:=d_{i+1}$ and degeneracy maps $s_i:=s_{i+1}$.
\end{defn}
 
\begin{observation}\label{simplicial path space}
Let $A_\bullet$ be a simplicial space and let $A_0$ denote the constant simplicial space. There are 
simplicial maps $\iota :A_0\rar A_\bullet$ and $\rho: PA_{\bu}\rar A_0$ defined on level $n$ via the maps 
$[n+1]\rar [0]$ and $[0]\hrar [n]$($0\mapsto 0$), respectively. $PA_{\bu}$ is simplicially homotopy 
equivalent to the constant simplicial space $A_0$; in particular, $|PA_\bu|\simeq A_0$. In addition, the face map $d_0:A_{n+1}\rar A_n$
defines a simplicial map $PA_\bu\rar A_\bu$.
\end{observation}
In addition, we will need the following result.
\begin{lem}\label{cartesian}
Let $\pi:A_\bu\rar B_\bu$ be a homotopy action. Then for each $n\geq 0$, the square 
\begin{equation*}
\xymatrix{A_{n+1}\ar[r]\ar[d] & |PA_\bu|\ar[d]\\A_n\ar[r] & |A_\bu|}
\end{equation*}
is homotopy cartesian.
\end{lem}
\begin{proof}
From the axioms of a homotopy action, there is a commutative square with horizontal maps homotopy equivalences
\begin{equation*}
\vcenter{
\xymatrix{A_{n+1}\ar[d]_{d_0}\ar[rrr]^{(d_0\cdots d_0)\times \pi_{n+1}} & & & 
A_0\times B_{n+1}\ar[d]^{1\times d_0} \\ 
A_n\ar[rrr]^{(d_0\cdots d_0)\times \pi_n} & & & A_0\times B_n.}}
\end{equation*}
Since $B_\bu$ is a reduced Segal space, it follows from \cite[1.6]{Seg} that for each $k\geq 0$, the square
\begin{equation}
\vcenter{
\xymatrix{B_{k+1}\ar[r]\ar[d]_{d_0} & |PB_\bu|\ar[d]\\ B_k\ar[r] & |B_\bu|}}
\end{equation}
is homotopy cartesian.

Thus, the homotopy fiber of $d_0:B_{n+1}\rar B_n$ is (canonically) equivalent to $B_1$. The homotopy fiber of 
$d_0:A_{n+1}\rar A_n$ is therefore homotopy equivalent to $B_1$, which is also the homotopy fiber of 
$|PA_\bu|\rar |A_\bu|$. 
It follows that the square
\begin{equation*}
\xymatrix{A_{n+1}\ar[r]\ar[d]_{d_0} & |PA_\bu|\ar[d]\\ A_n\ar[r] & |A_\bu|}
\end{equation*}
is homotopy cartesian.
\end{proof}

\begin{proof}[Proof of theorem \ref{homotopy action}]

\begin{enumerate}[\upshape (i)]
\item Given, without loss of generality, a fibration sequence \mbox{$X\!\!\rar\!X/G\!\rar\!\!BG$,} 
the map $hPow_\bu(X\rar X/G)\rar hPow_\bu(*\rar BG)$, obtained just as in theorem \ref{classification of fibrations} has 
$X$ as a homotopy fiber in each level. Since $|hPow_\bu(X\rar X/G)|\simeq X/G$ and 
$|hPow_\bu(*\!\rar\! BG)|\!\simeq \!BG$, the map 
\mbox{$|hPow_\bu(X\!\rar\! X/G)|\!\rar\! |hPow_\bu(*\!\rar\! BG)|$} 
is equivalent to $X/G\rar BG$.

\item Given a homotopy action $\pi: A_\bu\rar B_\bu$, $B_\bu$ is a reduced Segal space, 
and thus by proposition 1.6 in \cite{Seg}, for each $k\geq 0$, the square
\begin{equation}
\vcenter{
\xymatrix{B_{k+1}\ar[r]\ar[d]_{d_0} & |PB_\bu|\ar[d]\\ B_k\ar[r] & |B_\bu|}}
\end{equation}
is homotopy cartesian. By lemma \ref{cartesian}, the same holds for $A_\bu$, i.e.
for each $k\geq 0$, the square 
\begin{equation}
\vcenter{
\xymatrix{A_{k+1}\ar[r]\ar[d]_{d_0} & |PA_\bu|\ar[d]\\A_k\ar[r] & |A_\bu|}}
\end{equation}
is homotopy cartesian.
We construct a map $A_\bu \rar  hPow_\bu(|PA_\bu|\rar |A_\bu|)$ by 
induction on $n$. For $n\!=\!0$, the map 
$A_0 \!\rar \!|PA_\bu|$ is the realization of \mbox{$\iota:\! A_0 \rar PA_\bu$} defined in 
\ref{simplicial path space}. For $n=1$, consider the commutative square 
\begin{equation}
\vcenter{
\xymatrix{A_0\ar[r]\ar[d] & |A_\bu|\ar[d]\\ |PA_\bu|\ar[r] & |A_\bu|.}}
\end{equation} 
Since $(2)$ is homotopy cartesian for $k=0$, the map $A_1\rar A_0\times^h_{|A_\bu|} |PA_\bu|$ is a 
homotopy equivalence, and the map $A_1\rar hPow_1(|PA_\bu|\rar |A_\bu|)$ is obtained by composing the 
last map with $A_0 \times^h_{|A_\bu|} |PA_\bu|\rar |PA_\bu| \times^h_{|A_\bu|} |PA_\bu|$ induced by 
$(3)$. Let us define the map for $n+1$: the square $(2)$ with index $n$ is homotopy cartesian, 
and thus there is a homotopy equivalence $A_{n+1} \rar A_n\times^h_{|A_\bu|} |PA_\bu|$. 
Using the map $A_n\rar hPow_n(|PA_\bu|\rar |A_\bu|)$ that was defined, we get a natural homotopy equivalence $A_{n+1}\rar hPow_\bu(|PA_\bu|\rar |A_\bu|)$. 
It is clear from the construction that one gets a simplicial map 
$A_\bu \rar hPow_\bu(|PA_\bu|\rar |A_\bu|)$. Similarly, there is a simplicial map 
$B_\bu \rar hPow_\bu(|PB_\bu|\rar |B_\bu|)$. The zig-zag of commutative squares 
\begin{equation*}
\xymatrix{A_0\ar[r]^\simeq\ar[d] &|PA_\bu|\ar[r]\ar[d] & |PB_\bu|\ar[d] &  B_0\ar[l]_\simeq\ar[d]\\ 
|A_\bu|\ar[r]^\simeq & |A_\bu|\ar[r] & |B_\bu| & |B_\bu|\ar[l]_\simeq}
\end{equation*}
induces  a zig-zag of commutative simplicial squares
\begin{equation*}
\vcenter{
\xymatrix{A_\bu\ar[r]^(0.2){\simeq}\ar[d] & hPow_\bu(|PA_\bu|\rar |A_\bu|)\ar[d] & 
hPow_\bu(A_0 \rar |A_\bu|)\ar[l]_(0.45){\simeq}\ar[d]\\ 
B_\bu\ar[r]^(.2){\simeq} & hPow_\bu(|PB_\bu| \rar |B_\bu|) & hPow_\bu(B_0\rar |B_\bu|).\ar[l]_(0.45){\simeq} }}
\end{equation*}
Note that by proposition \ref{Bar as IFP}, there is also a square
\begin{equation*}
\xymatrix{hPow_\bu(A_0 \rar |A_\bu|)\ar[d] & Bar_\bu(X,G)\ar[d]\ar[l]\\ 
hPow_\bu(B_0\rar |B_\bu|) & B_\bu(G)\ar[l]}
\end{equation*}
for a topological group $G$ with $BG\simeq |B_\bu|$.

\end{enumerate}
\end{proof}

\section{An invariant characterization of normality}
Theorem \ref{[FS]} characterizes homotopy normal maps of discrete groups in terms of a simplicial 
group, equivariantly equivalent to the bar construction. By analogy, the mere fact that the homotopy bar 
construction $Bar_\bullet(\Om Y,\Omega X)$ is simplicially equivalent to a simplicial loop space 
$\Gamma_\bullet$ with $\Gamma_0\simeq\Omega Y$, is a necessary but not sufficient condition for a loop map 
$\Omega f:\Omega X\rar \Omega Y$ to be homotopy normal.

In both simplicial spaces $Bar_\bullet(\Omega Y,\Omega X)
$ and $Q_\bullet$
(see notation \ref{standard notation}), the map \linebreak
$s_{n-1}\cdots s_0$ is a loop map, 
and therefore it induces a homotopy action of $\Omega Y$ on 
 $Q_n$ and $Bar_n(\Omega Y,\Omega X)$ (see corollary \ref{induced action}). 

We begin with the following.
\begin{prop}\label{equivalence of actions}
Let $\Omega f:\Omega X\rar \Omega Y$ be a homotopy normal map and $Q_\bu$ its corresponding simplicial loop space. For each $n$, the homotopy actions 
induced by the loop maps $Q_0\rar Q_n$ and 
$\Omega Y\rar Bar_n(\Omega Y,\Omega X)$ are 
equivalent via the map $\eta:Q_\bullet\rar Bar_\bullet(\Omega Y,\Omega X)$, defined in 
\ref{standard notation}. 
\end{prop}
\begin{proof}
We do only the case $n=1$ since other cases are similar. 
Write \linebreak
$\sigma:=s_0:Q_0\rar Q_1$ and 
$s:=s_0:Bar_0(\Omega Y,\Omega X)\rar Bar_1(\Omega Y,\Omega X)$. 
The simplicial equivalence $\eta:Q_\bullet\rar Bar_\bullet(\Omega Y,\Omega X)$
induces a commutative square with vertical arrows being homotopy 
equivalences, and with the left vertical arrow being a loop map 
\begin{equation*}
\vcenter{
\xymatrix{Q_0\ar[d]_{\eta_0}\ar[r]^{\sigma} & Q_1\ar[d]^{\eta_1}\\ 
\Omega Y\ar[r]^{s} &  \Omega Y\times\Omega X.}}
\end{equation*}
Finding the dashed arrow 
\begin{equation*}
\vcenter{
\xymatrix{Q_1\ar[r]^\gamma\ar[d]_{\eta_1} & Q_1//Q_0\ar@{-->}[d]^{d_1}\\ 
\Omega Y\times\Omega X\ar[r]_c & \Omega X}}
\end{equation*} 
will end the proof because the first and second homotopy actions are built out of homotopy powers of 
$\gamma$ and $c$, respectively. Both $\sigma$ and $s$ have 
(spaces of the homotopy type of) loop spaces as their homotopy fiber, and the Puppe-Nomura sequence will 
provide the dashed arrow, once we show that the equivalence between the homotopy fibers 
$F:=hfib(\sigma)\rar hfib(s)\simeq\Omega^2 X$ is a loop map. To prove the last statement we use 
the path-space to model the homotopy fiber. On the one hand, we have the pullback square 
\begin{equation*}
\vcenter{
\xymatrix{{\Omega}^2 X\ar[r]\ar[d] & P( \Omega Y\times\Omega X)\ar[d] \\ 
\Omega Y\ar[r]^s & \Omega Y\times \Omega X,}}
\end{equation*}
and on the other hand, in the pullback square 
\begin{equation*}
\vcenter{
\xymatrix{F\ar[r]\ar[d] & P(Q_1)\ar[d] \\ Q_0\ar[r]^\sigma & Q_1,}}
\end{equation*}
all maps are of the homotopy type of loop maps. The map $F\rar \Om^2 X$ is the universal map to the pullback $\Om^2 X$, obtained from the diagram
$$\xymatrix{F\ar@{.>}[dr]\ar@/^1.5pc/[drr]\ar@/_1.5pc/[ddr] & & \\ & \Om^2 X\ar[r]\ar[d] & P(\Om Y\times\Om X)\ar[d]\\ & \Om Y\ar[r] & \Om Y\times\Om X }$$
where the curved maps are
$F\rar Q _0\rar \Omega Y $ and $F\rar P(Q_1)\rar P(\Omega Y\times\Omega X)$; these maps are (of the homotopy type of) loop maps, and thus the map they induce 
$F\rar \Omega^2 X$ is itself (of the homotopy type of) a loop map.
\end{proof}

As we have just seen, the loop maps 
$s_{n-1}\cdots s_0: Q_0\rar Q_n$ ($n=0$ understood as the identity map) 
induce homotopy actions of $Q_0$ on $Q_n$. We can pack all the maps into one simplicial map 
$Q_0\rar Q_\bullet$, which will then induce a simplicial object in the category of homotopy 
actions. Recalling definition \ref{homotopy action on a simplicial space}, 
this is a homotopy action of $Q_0$ on $Q_\bullet$. Similarly, one has a homotopy action of 
$\Omega Y$ on $Bar_\bullet(\Omega Y,\Omega X)$ and the loop 
space equivalence  $Q_0\simeq \Omega Y$ makes the 
first homotopy action into one of $\Omega Y$ on $Q_\bullet$
(see \ref{induced action}). Note that any simplicial loop space $\Gamma_\bullet$ with $\Gamma_0\simeq \Om Y$ could play the role of $Q_\bullet$ in defining these homotopy actions. 

Given a loop map $\Om f:\Omega X\rar \Omega Y$ and a simplicial loop space $\Gamma_\bullet$ with $\Gamma_0\simeq \Om Y$, we call the actions above the \emph{canonical homotopy actions} of $\Omega Y$ on $\Gamma_\bullet$ and $Bar_\bullet(\Omega Y,\Omega X)$. 
The additional condition for a characterization of normality is that the two are equivalent.

We can now restate and prove theorem A.
\begin{thm}\label{characterization of homotopy normality}
A loop map $\Om f:\Omega X\rar \Omega Y$ is homotopy normal 
if and only if there exist a simplicial loop space $\Gamma_\bullet$ with $\Gamma_0 \simeq \Omega Y$ (as loop spaces), 
and such that the canonical homotopy actions of $\Omega Y$ on $\Gamma_\bullet$ and on
$Bar_\bullet(\Omega Y,\Omega X)$ (as above) are weakly equivalent.
\end{thm}
\begin{rem}
The weak equivalence of homotopy actions above implies, in particular, the equivalence of simplicial spaces $Bar_\bu(\Om Y,\Om X)$ and $\Gamma_\bu$.
\end{rem}
\begin{proof}
Assume $\Omega f$ is homotopy normal. 
We have a commutative square of simplicial spaces 
\begin{equation*}
\vcenter{
\xymatrix{\Omega Y\ar[r]^{\sigma}\ar[d]^{1} & Q_\bullet\ar[d]^{\varphi}\ar[r] & 
Q_\bullet//Q_0\ar@{-->}[d]^d\\ 
\Omega Y\ar[r]^<(0.2){s} & Bar_\bullet(\Omega Y,\Omega X)\ar[r] & 
Bar_\bullet(\Omega Y,\Omega X)//\Omega Y}}
\end{equation*}
with $\varphi$ the simplicial equivalence of theorem \ref{Loop structure}; 
the dashed arrow $d$ with $d_1$ (of proposition \ref{equivalence of actions}) as its first component, 
and the analogous $d_n$ as its $n$-th component. 
This gives the desired equivalence of the canonical actions.\\
Conversely, if we have a zig-zag of equivalent homotopy actions (see \ref{homotopy action on a simplicial space}), then taking the homotopy quotient of each homotopy action, we get a zig-zag of simplicial spaces
\begin{equation*}
\vcenter{
\xymatrix{\Gamma_0\ar[r]\ar[d]^\simeq & \Gamma_\bu\ar[r]\ar[d]^\simeq & \Gamma_\bu//\Gamma_0\ar[d]^\simeq\\ 
\vdots\ar[r] & \vdots\ar[r] & \vdots\\
\Om Y\ar[r]\ar[u]_\simeq & Bar_\bu(\Omega Y,\Omega X)\ar[r]^<(0.2)q\ar[u]_\simeq & Bar_\bu(\Om Y,\Om X)//\Om Y.\ar[u]_\simeq}}
\end{equation*}
The map $q$ in the bottom row is in fact $\pi:Bar_\bu(\Om Y,\Om X)\rar Bar_\bu(*,\Om X)$, 
and upon realization we have a zig-zag of equivalent principal fibrations 
\begin{equation*}
\vcenter{
\xymatrix{\Gamma_0\ar[r]\ar[d]^\simeq & |\Gamma_\bullet|\ar[r]\ar[d]^\simeq & |\Gamma_\bu//\Gamma_0|\ar[d]^\simeq\\
\vdots\ar[r]& \vdots\ar[r] & \vdots\\
\Omega Y\ar[r]\ar[u]_\simeq & \Omega Y//\Omega X\ar[r]\ar[u]_\simeq & X.\ar[u]_\simeq }}
\end{equation*}
The operation of taking loops commutes with that of realization, and hence 
$|\Gamma_\bullet|\simeq \Omega W$ for some connected space $W$. 
The map $\Gamma_0\rar |\Gamma_\bullet|$ is the realization of a simplicial loop map 
$\Gamma_0\rar \Gamma_\bullet$, hence a loop map itself, 
and delooping it gives the desired extension $Y\drar W$.
\end{proof}
As an application of theorem \ref{characterization of homotopy normality} we will show that homotopy normal maps are preserved by HM functors.

Let $A_\bullet \rar B_\bullet$ be a homotopy action. From proposition \ref{correspondence} (b), 
it follows that there is a homotopy fibration sequence 
$A_0\overset{\sigma}{\rar} |A_\bullet|\rar |B_\bullet|$, where $\sigma$ is the realization of the 
simplicial map $A_0\rar A_\bullet$ that has as $n$-th component the map $s_{n-1}\cdots s_0$. Since 
$B_\bullet$ is a reduced Segal space, $\Omega |B_\bullet|\simeq B_1$. We denote by 
$\psi:B_1\rar A_0$ the canonical map from the homotopy fiber of 
$\sigma:A_0\rar |A_\bullet|$ to $A_0$ and endow $A_0$ with a base-point via $\psi$. 
Denote by $i:B_1\rar A_0\times B_1$ the natural inclusion.
We shall need the following technical lemma.
\begin{lem}\label{Tech}
For any choice of homotopy inverse $e:A_0\times B_1 \rar A_1$ for \linebreak
$d_1\times \pi_1:A_1\rar  A_0\times B_1$, the composite 
$B_1\overset{i}{\rar}  A_0\times B_1\overset{e}{\rar} A_1\overset{d_0}{\rar} A_0$ 
is homotopic to $\psi$.
 \end{lem}
\begin{proof}
The square 
\begin{equation}
\vcenter{
\xymatrix{A_1\ar[r]^{d_1}\ar[d]_{d_0} & A_0\ar[d] \\ A_0\ar[r] & |A_\bullet|}}\tag{$\ast$}
\end{equation}
is homotopy commutative. We thus obtain a homotopy commutative diagram of solid arrows
\begin{equation*}
\vcenter{
\xymatrix{B_1\ar[r]^{i} & A_0\times B_1\ar[r]^{pr}\ar@/^/@{-->}[d]^e & A_0
\\ B_1\ar[r]\ar[d]^{\simeq}_{c_1}\ar[u]_{\simeq}^{c_2} & A_1\ar[r]^{d_1}\ar[d]^{d_0}\ar[u]^{d_1\times \pi_1} & A_0\ar[d]^{\sigma}\ar[u]^{1}
\\ B_1\ar[r]^{\psi} & A_0\ar[r]^{\sigma} & |A_{\bullet}|, }}
\end{equation*}
where the map $B_1\rar A_1$ is the canonical map from the homotopy fiber, the map $c_1$ is the comparison map between the homotopy fibers of 
$d_1$ and $\sigma$, which is a homotopy equivalence, and the map $c_2$ is the comparison map between the homotopy fibers of $d_1$ and $pr$, 
which is again a homotopy equivalence.  The lemma now follows from inverting $c_2$.
 \end{proof}

\begin{thm}\label{invariance of bar}
Let $\Om f:\Om X\rar \Om Y$ be a loop map and \mbox{$L\!:\!Top\!\rightarrow\! Top$} an HM functor. Then $LBar_\bu(\Om Y,\Om X)\rar LBar_\bu(*,\Om X)$ is weakly equivalent to $Bar_\bu(L\Om Y,L\Om X)\rar Bar_\bu(*,L\Om X)$ where the latter is induced from $L\Om f$. 
\end{thm}
\begin{proof}

$L\Omega Y\rar |LBar_\bullet(\Omega Y,\Omega X)|\rar |LBar_\bullet(*,\Omega X)|$ 
is a homotopy fibration sequence (being the realization of a simplicial fibration sequence), 
and since $| LBar_\bullet(*,\Omega X)|\simeq B(L\Omega X)$ ($LBar_\bullet(*,\Omega X)$ 
is a reduced Segal space for $L\Omega X$), there is a map 
$\varphi: L\Omega X\rar L\Omega Y$, which is the map from the homotopy
fiber of $L\Omega Y\rar |LBar_\bullet(\Omega Y,\Omega X)|$ to $L\Om Y$.

Abbreviate $A_\bullet:=Bar_\bullet(\Omega Y,\Omega X)$ and $B_\bullet:=Bar_\bu(*,\Omega X)$.  
If $e:A_0\times B_1\rar A_1$ is a homotopy inverse to $d_1\times \pi_1$, 
then $Le$ is a homotopy inverse for $L(d_1\times \pi_1)$, which is
equivalent to $L(d_1)\times L(\pi_1)$. By lemma \ref{Tech}, 
$\Omega f$ is homotopic to the composite 
$B_1\overset{i}{\rar} A_0\times B_1 \overset{e}{\rar} A_1 \overset{d_0}{\rar} A_0$, 
and so $L\Omega f$ is homotopic to the composition $Ld_0\circ Le\circ Li$. 
The last composite is homotopic to the composite 
$LB_1\hrar  LA_0\times LB_1 \overset{Le\circ w}{\lrar} LA_1\overset{Ld_0}{\rar} A_0$ (where $w$ is some homotopy inverse for $L(A_0\times B_1)\rar  LA_0\times LB_1$),
which is in turn homotopic to $\varphi$ by lemma \ref{Tech} ($Le\circ w$ is a homotopy inverse for $Ld_1\times L\pi_1$). 
It follows that $L\Omega f$ is equivalent to $\varphi$. Thus, the map $L\Om Y\rar |LBar_\bu(\Om Y,\Om X)|$ is equivalent to $L\Om Y\rar L\Om Y//L\Om X$ and using proposition \ref{Bar as IFP} and theorem \ref{homotopy action}, we deduce that $LBar_\bullet(\Om Y,\Om X)\rar LBar(*,\Om X)$ is weakly equivalent to $Bar_\bullet(L \Om Y,L\Om X)\rar Bar_\bullet(*,L\Om X)$.
\end{proof}
Let us rephrase theorem \ref{invariance of bar}. Given a loop map $\Omega f$ and an HM functor $L$, there are two homotopy actions: the first is given by applying $L$ to the homotopy action induced by $\Om f$, and the second is the homotopy action induced from $L\Om f$. The theorem then says that the two are weakly equivalent. We note that if we are given a homotopy action of a loop space on a simplicial space, in which the homotopy actions in each level are induced by loop maps, an analogous statement holds.
 
Using the machinery of reduced Segal spaces, one can easily see that applying an HM functor to a 
simplicial loop space in every level yields a simplicial space \emph{simplicially equivalent} to a 
simplicial loop space.  
 
Thus, we now know all the ingredients used in theorem \ref{characterization of homotopy normality} are invariant under HM functors and we deduce theorem B (which we restate for convenience).
\begin{thm}[Theorem B]\label{invariance}
Let $\Om f:\Omega X\rar \Omega Y$ be a homotopy normal map. 
If \mbox{$L\!:\!Top\!\rightarrow\! Top$} is an HM functor, then $\L (\Om f): L\Omega X\rar L\Omega Y$ 
is a homotopy normal map.
\end{thm}

Let us demonstrate a use of theorem \ref{invariance} by applying it to prove theorem C (which we restate).
\begin{thm}\label{tc}
Let $p:E\rar B$ be a principal fibration with $B$ connected, 
$f:X\rar Y$ a map of pointed connected spaces and $L_{\Sigma f}$ 
the localization with respect to its suspension. 
Then $L_{\Sigma f}E\rar L_{\Sigma f}B$ is equivalent to a principal fibration.
\end{thm}
\begin{rem}
Note that if $G$ is the structure group of $E\rar B$, $L_{\Sigma f}G$ 
need not be the structure group of $L_{\Sigma f}E\rar L_{\Sigma f}B$.
\end{rem}
\begin{proof}[Proof of Theorem \ref{tc}]
Note that $\Omega E\rar \Omega B$ is homotopy normal. Hence, \linebreak 
$L_{f}\Omega E\rar L_{f}\Omega B$ is homotopy normal. 
Since for any pointed space $A$ there is a natural equivalence 
$L_{f}\Omega A\simeq \Omega L_{\Sigma f} A$, we get that 
$\Omega L_{\Sigma f} E\rar \Omega L_{\Sigma f} B$ is homotopy normal and thus 
$L_{\Sigma f}E\rar L_{\Sigma f}B$ is a homotopy principal fibration.
\end{proof}

\section{Higher normality}
As mentioned in example \ref{double loop map}, any double loop map with simply-connected underlying 
spaces is automatically homotopy normal. 
However, in the case of a double loop map, it is more natural 
to ask when the homotopy quotient admits a natural double loop space structure. 
 \begin{defn}
A $0$-homotopy normal map is a pointed map which admits a structure 
of a (homotopy) principal fibration of connected spaces. For $k\geq 1$, 
call a $k$-fold loop map ${\Omega}^k f:{\Omega}^k X\rar {\Omega}^k Y$ \emph{k-homotopy normal} if $f$ is $0$-homotopy normal. 
 \end{defn}

Thus, if a k-fold loop map ${\Om}^k f$ is k-homotopy normal, the homotopy quotient ${\Omega}^k Y//{\Omega}^k X$ (which is always a (k-1)-fold loop space) admits a structure of a k-fold loop space in a natural way. 

\begin{rem}
One may wonder about the definition of `$\infty$-homotopy normality'. 
However, any infinite loop map $X\rar Y$ induces a principal fibration sequence 
of infinite loop spaces $X\rar Y\rar Y//X$. Thus any infinite loop map is `$\infty$-normal' in the naive sense. 
This is a reflection of the fact that any inclusion map of abelian (topological) 
groups is the inclusion of a normal subgroup. 

\end{rem}

We begin with an extension of theorem A.
\begin{thm}\label{characterization of k-normality}
A k-fold loop map ${\Omega^k f}:\Omega^k X{\lrar} \Omega^k Y$ is 
$k$-homotopy normal if and only if there exists a $k$-fold simplicial loop space 
$\Gamma_\bullet$ with $\Gamma_0\simeq\Omega^k Y$, 
and such that the canonical homotopy actions of $\Omega^k Y$ on 
$Bar_\bullet(\Omega^k Y,\Omega^k X)$ and $\Gamma_\bullet$ are naturally equivalent.  
\end{thm}
\begin{proof}
This is analogous to the proof of theorem \ref{characterization of homotopy normality}. 
If $\Omega^k f$ is $k$-homotopy normal, then $\Omega f$ is homotopy normal, and 
looping down its extension $Y\rar W$ $k$ times gives a $k$-fold loop map equivalent to 
$\Omega^k Y\rar \Omega^k Y//\Omega^k X$. Taking the (homotopy) power of that map gives 
the desired $k$-fold loop space. Conversely, such a $k$-fold loop space 
gives a (homotopy) principal fibration sequence of $k$-fold loop spaces 
$\Omega^k X\rar \Omega^k Y\rar |\Gamma_\bullet|$, equivalent to the Borel construction, 
providing the $k$-homotopy normality required.
\end{proof}
We wish to use the same methods as before to prove invariance of $k$-homotopy normal 
maps under $HM$ functors. For that, we need to know that $k$-fold loop spaces 
are invariant under these functors. A slight generalization of 
reduced Segal spaces is the tool needed. 
\begin{defn}
Let $k$ be a positive integer. A \emph{k-simplicial space} is a functor 
$$
\Delta^{op} \times\cdots \times\Delta^{op}\rar Top \;\;(k\;\;times)
$$
\end{defn}

The following is taken from \cite{BFSV}

\begin{defn}
A  k-simplicial space $X_{\bullet\cdots\bullet}$ is called a \emph{reduced Segal k-space} if:
\begin{enumerate}
\item
  $X_{0,\ldots,0}\simeq *$;
\item
the Segal maps induce homotopy equivalences 
$X_{p_1,\ldots,p_k}\overset{\simeq}{\rar} (X_{1,\ldots,1})^{p_1\cdots p_k}$;
\item the monoid $\pi_0(X_{1,\ldots,1})$ admits inverses (i.e. is a group).
\end{enumerate}
\end{defn} 

Building on Segal's delooping machine, the characterization of $k$-fold loop spaces takes the following form.
\begin{thm}
A space $X$ is of the homotopy type of a $k$-fold loop space if and only if there exist a reduced Segal k-space $X_{\bullet,\ldots,\bullet}$ with $X_{1,\ldots,1}\simeq X$.
\end{thm} 

\begin{cor}
Homotopy monoidal endofunctors of spaces preserve $k$-fold loop spaces.
\end{cor}

Using exactly the same arguments of theorem \ref{invariance}, theorem \ref{characterization of k-normality} implies that $L$ preserves higher homotopy normality.
\begin{thm}
If ${\Omega}^k f:{\Omega}^k X\rar {\Omega}^k Y$ is $k$-homotopy normal and 
$L:Top\rar Top$ an HM functor, then $L(\Omega^k f)$ is $k$-homotopy normal.
\end{thm}

\end{document}